\documentclass[11pt,a4paper]{amsart}
\pdfoutput=1

\usepackage[utf8]{inputenc}

\usepackage{hyperref,graphicx,amssymb}

\DeclareMathOperator{\Gal}{Gal}

\newcommand{\CC}{\mathbb{C}}
\newcommand{\PP}{\mathbb{P}}
\newcommand{\HH}{\mathbb{H}}

\begin{document}

\title{Belyi map for the sporadic group $J_1$}

\author{Dominik Barth}
\author{Andreas Wenz}

\address{Institute of Mathematics\\ University of Würzburg \\ Emil-Fischer-Straße 30 \\ 97074 Würzburg, Germany}
\email{dominik.barth@mathematik.uni-wuerzburg.de}
\email{andreas.wenz@mathematik.uni-wuerzburg.de}



\begin{abstract}
We compute the genus $0$ Belyi map for the sporadic Janko group $J_1$ of degree $266$ and describe the applied method. This yields explicit polynomials having $J_1$ as a Galois group over $K(t)$, $[K:\mathbb{Q}] = 7$.
\end{abstract}
\maketitle

\section{Introduction}

Our main goal of this paper is to present a method to compute high degree genus $0$ Belyi maps having prescribed monodromy. Using this technique we found explicit polynomials having the sporadic Higman-Sims group as Galois groups over $\mathbb{Q}(t)$, see \cite{Barth2016}, and more generally all Belyi maps of genus $0$ having almost simple, primitive monodromy groups (not $A_n$, $S_n$) generated by rigid rational triples of degree between $50$ and $250$, see \cite{Barth2017}.
In the following we will apply this method to compute a Belyi map of degree $266$ for the permutation triple $(x,y,z:=(xy)^{-1}) \in S_{266}^3$, given in the ancillary data. It is of type
\begin{center}
\renewcommand{\arraystretch}{1.3}
\begin{tabular}{c|c|c|c}
& $x$ & $y$ & $z = (xy)^{-1}$ \\  \hline
cycle structure & $7^{38}$ & $2^{128}.1^{10}$ & $3^{87}.1^5$ \\
\end{tabular}
\end{center}
and has the following properties:
\begin{itemize}
\item $x,y$ generate the sporadic Janko group $J_1$.
\item $(x,y,z)$ is of genus $0$.
\item The permutations $x,y,z$ each lie in rational conjugacy classes.
\item $(x,y,z)$ has passport size $7$.
\end{itemize}

Alternative techniques for computing Belyi maps can be found in \cite{Sijsling2014}, \cite{Roberts2016}, \cite{Vidunas2016}, \cite{Monien2014} and \cite{Klug2014}. Recently, Monien \cite{Monien2017} demonstrated another powerful method by computing a Belyi map for $J_2$ of degree $100$.

\section{Method of Computation}

Let $a:= \text{ord}(x)$, $b:= \text{ord}(y)$, $c:= \text{ord}(z)$ and
$$
\Delta:=\left< \delta_a,\delta_b,\delta_c \,|\, \delta_a^a= \delta_b^b=\delta_c^c = \delta_a\delta_b\delta_c = 1 \right>.
$$
Note that $(x,y,z)$ is hyperbolic since $1/a+1/b+1/c < 1$.
We now consider the embedding $\Delta \hookrightarrow \text{PSL}_2(\mathbb{R})$ described in \cite[Proposition 2.5]{Klug2014}, where $\delta_a$ (resp. $\delta_b$) is mapped to a hyperbolic rotation around $i$ (resp. $\mu i$ for some $\mu >1$) of angle $\pi/a$ (resp. $\pi/b$).
Thus $\Delta$ acts on the upper half-plane $\mathbb{H}$ via the natural action of $\text{PSL}_2(\mathbb{R})$ on $\mathbb{H}$ and its fundamental domain is the hyperbolic kite with vertices $i,\gamma,\mu i,-\gamma$ for some $\gamma \in \mathbb{H}$. Furthermore, let $\varphi$ denote the homomorphism from $\Delta$ onto $G:=\langle x,y \rangle$ such that $\delta_a \mapsto x$ and $\delta_b\mapsto y$.
With the notation $\Gamma := \varphi^{-1}(G_1)<\Delta$ we define
$$\Phi:\mathbb{H}/\Gamma \to \mathbb{H}/\Delta, \; z \text{ mod } \Gamma \mapsto z \text{ mod } \Delta. \footnote{For a topological space $X$ and a group $G$ acting on $X$ let $X/G$ denote the corresponding orbit space.}$$
This is a Belyi map, i.e. a three branch point covering, of degree $d$ ramified over $i, \mu i$ and $\gamma$ and its monodromy group is isomorphic to $G$, see for example \cite{Klug2014} for more details.

One can now restrict $\Phi$ to a connected fundamental domain $D \subset \HH$ by introducing an equivalence relation $\sim$ on $\partial D$ induced by the quotient structure of $\mathbb{H}/\Gamma$. After identifying $\mathbb{H}/\Delta$ with $\mathbb{P}^1(\mathbb{C})$ such that $i\mapsto 0$, $\mu i \mapsto 1$ and $\gamma \mapsto \infty$ it suffices to study the following induced Belyi map
$$
F: D/_\sim\, \to \mathbb{P}^1(\mathbb{C}).
$$
In the following we describe how to get a Belyi map of type $\mathbb{P}^1(\mathbb{C}) \to \mathbb{P}^1(\mathbb{C})$, this is possible because our given permutation triple is of genus $0$, thus $D/_\sim \cong \PP^1(\CC)$.
The Schwarz-Christoffel Toolbox \cite{Driscoll1996} for MATLAB \cite{MATLAB} gives us an approximation of a conformal map $$
h_1:D^\circ \to \mathbb{H}.$$
We can extend $h_1$ on $\overline{D}$ such that $h_1(\partial D) = \partial \mathbb{H}$ and $h_1$ is well defined on $D/_\sim$. Now $\sim$ induces via $h_1$ a new equivalence relation $\approx$ on $\partial \mathbb{H}$.
In order to glue the corresponding edges on $\partial \mathbb{H}$ we apply the slit algorithm (see  \cite{Marshall2007} and \cite{Barnes2014}). We therefore find a tree $T$ in $\mathbb{P}^1(\mathbb{C})$ and a conformal map 
$$
h_2:\mathbb{H} \to \mathbb{P}^1(\mathbb{C}) \setminus T $$
such that
$h_2$ can be extended on $\overline{\mathbb{H}}$, $h_2(\partial\mathbb{H}) = T$ and $h_2$ is well defined on $\mathbb{H}/_\approx$.
This leads to the Belyi map
$$F\circ h_1^{-1} \circ h_2^{-1} : \mathbb{P}^1(\mathbb{C}) \to \mathbb{P}^1(\mathbb{C})$$
having the prescribed ramification data $(x,y,z)$ over $0$, $1$ and $\infty$.

Practically, the above method yields approximate preimages of $0$, $1$ and $\infty$. Using Newton's method we compute this Belyi map with sufficiently high precision, allowing us to recognize its coefficients as algebraic numbers using the LLL algorithm \cite{Lenstra1982} implementation provided by Magma~\cite{Magma}.

\section{Computed Results and Verification}

The computed Belyi map $$f = \frac{p}{q} = 1 + \frac{r}{q},$$ defined over the number field 
$$
K = \mathbb{Q}(\alpha) \text{ with } \alpha^7 - \alpha^6 - 2\alpha^4 - \alpha^3 + 2\alpha^2 + 2\alpha + 2=0
$$	
can be found in the ancillary file and its reduction modulo a prime ideal of norm 269 is given in the appendix. 

With the Riemann-Hurwitz genus formula and the factorizations of $p,q$ and $r$ one can easily verify that $f:\mathbb{P}^1(\CC) \rightarrow \mathbb{P}^1(\CC)$ is indeed a Belyi map, i.e. a three branch point covering, ramified over $0$, $1$ and $\infty$.

\begin{proof}[Verification of monodromy]

Let $\mathfrak{p}=(5,2+\alpha)\mathcal{O}_K$ be the unique prime ideal in the ring of integers $\mathcal{O}_K$ of $K$ with norm $5$. Then all coefficients of $p(X)-tq(X) \in K(t)[X]$ lie in the localization $R:=\mathcal{O}_K[t]_{\mathfrak{p}\mathcal{O}_K[t]}$ of $\mathcal{O}_K[t]$ at the prime ideal $\mathfrak{p}\mathcal{O}_K[t]$. Denote by $\bar{p}$ and $\bar{q}$ the images under the canonical homomorphism $R[X] \rightarrow (R/\mathfrak{p}R)[X] \cong \mathbb{F}_5(t)[X]$.
The decomposition algorithm, found in \cite{Ayad2008}, implemented in Magma \cite{Magma}, yields that $\bar{p}/\bar{q} \in \mathbb{F}_{5}(X)$ is indecomposable, thus $A_{\mathbb{F}_5} := \Gal(\bar{p}(X)-t\bar{q}(X) \in \mathbb{F}_5(t)[X])$ is primitive by Ritt's Theorem.
Since $K(t)$ is the quotient field of $R$ we find $A_{\mathbb{F}_5} \leq A_K := \Gal(p(X)-tq(X) \in K(t)[X])$ using Dedekind reduction, see \cite[VII, Theorem 2.9]{Lang2002}. It follows that $A_K$ is primitive.

Furthermore, note that $(p(t)q(X) - q(t)p(X))/(X-t) \in K(t)[X]$ is reducible\footnote{The ancillary data contains a  factor of $p(t)q(X) - q(t)p(X)$ of degree $11$. This divisor was computed by using sufficiently enough specializations in $t$ and factorizing the resulting polynomials giving rise to an interpolation.}, thus $A_K$ is not $2$-transitive. Since $J_1$, $A_{266}$, and $S_{266}$ are the only primitive groups of degree 266 we find $A_K  = J_1$.
It is well-known that the geometric monodromy group $G := \Gal(p(X)-tq(X) \in \CC(t)[X])$  is normal in $A_K$.
Since $J_1$ is simple, we have $G = A_K = J_1$.
\end{proof}

\section*{Acknowledgements}

\noindent
We would like to thank Joachim König and Peter Müller for some helpful suggestions.

\appendix
\section*{Appendix: Computed Data}

In this section we will present the computed Belyi map $$ f = \frac{p}{q} = 1+ \frac{r}{q}$$ in the following way: We take the prime ideal $\mathfrak{p}=(269,207+\alpha)\mathcal{O}_K$ of norm $269$ and reduce the coefficients of $p,q \in K[X]$ by $\mathfrak{p}$ to obtain polynomials $\bar{p}$ and $\bar{q}$ lying in $(\mathcal{O}_K/\mathfrak{p})[X] \cong \mathbb{F}_{269} [X]$. The results are
\begin{align*}
\bar{p}(X) =\; &(X + 224)^7 \cdot\\&
(X^2 + 201X + 17)^7 \cdot \\&
(X^2 + 225X + 175)^7 \cdot \\&
(X^3 + 146X + 146)^7 \cdot \\&
(X^3 + 7X^2 + 88X + 26)^7 \cdot \\&
(X^3 + 176X^2 + 124X + 227)^7 \cdot \\&
(X^6 + 77X^5 + 257X^4 + 141X^3 + 224X^2 + 252X + 140)^7 \cdot \\&
(X^6 + 105X^5 + 148X^4 + 149X^3 + 231X^2 + 132X + 250)^7 \cdot \\&
(X^6 + 116X^5 + 61X^4 + 230X^3 + 119X^2 + 114X + 201)^7\cdot \\&
(X^6 + 257X^5 + 46X^4 + 100X^3 + 79X^2 + 188X + 127)^7
\end{align*}
and
\begin{align*}
\bar{q}(X) =\;  &(X^2 + 182X + 28)^3 \cdot\\&
(X^5 + 107X^4 + 235X^3 + 37X^2 + 91X + 15) \cdot \\&
(X^5 + 145X^4 + 46X^3 + 241X^2 + 163X + 209)^3 \cdot \\& 
(X^{10} + 29X^9 + 90X^8 + 199X^7 + 20X^6 + 220X^5 \\& + 18X^4 + 104X^3 + 212X^2 + 168X + 126)^3 \cdot \\&
(X^{10} + 41X^9 + 63X^8 + 219X^7 + 177X^6 + 124X^5 \\& + 134X^4 + 246X^3 + 206X^2 + 41X + 18)^3 \cdot \\&
(X^{10} + 55X^9 + 247X^8 + 253X^7 + 161X^6 + 49X^5 \\& + 242X^4 + 113X^3 + 235X^2 + 212X + 169)^3 \cdot \\&
(X^{10} + 67X^9 + 51X^8 + 12X^7 + 135X^6 + 116X^5 \\& + 172X^4 + 265X^3 + 239X^2 + 45X + 119)^3 \cdot \\&
(X^{10} + 76X^9 + 193X^8 + 151X^7 + 120X^6 + 172X^5 \\& + 192X^4 + 256X^3 + 29X^2 + 68X + 24)^3 \cdot \\&
(X^{10} + 98X^9 + 229X^8 + 244X^7 + 142X^6 + 180X^5 \\& + 65X^4 + 188X^2 + 23X + 227)^3 \cdot \\&
(X^{10} + 164X^9 + 124X^8 + 10X^7 + 268X^6 + 169X^5 \\& + 204X^4 + 222X^3 + 106X^2 + 30X + 200)^3 \cdot \\&
(X^{10} + 194X^9 + 187X^8 + 100X^7 + 204X^6 + 145X^5 \\& + 224X^4 + 67X^3 + 105X^2 + 203X + 5)^3 \cdot
\end{align*}

\end{document}